\documentclass[12pt,dvips]{amsart}
\usepackage{amsfonts, amssymb, latexsym, epsfig}
\usepackage{euler, amsfonts, amssymb, latexsym, epsfig}
\usepackage{euler, amsfonts, amssymb, latexsym, epsfig,epic}

\newtheorem{Theorem}{Theorem} 

\newtheorem{Proposition}{Proposition} 
\newtheorem{Lemma}{Lemma}

\newtheorem{Corollary}{Corollary}
\newtheorem*{Corollary*}{Corollary}
 
\newtheorem*{Theorem*}{Theorem}

\theoremstyle{remark}
\newtheorem{Example}{Example}

\newcommand\X{X}

\newcommand{\plainstack}[2]{\genfrac{}{}{0pt}{}{#1}{#2}}

\theoremstyle{plain}

\newtheorem{Problem}{Problem}
\newtheorem{Conjecture}{Conjecture}               

\font\co=lcircle10

\def\jr{\smash{\raise2pt\hbox{\co \rlap{\rlap{\char'005} \char'007}}
               \raise6pt\hbox{\rlap{\vrule height5pt}}
               \raise2pt\hbox{\rlap{\hskip4pt \vrule height0.4pt depth0pt
                width5.7pt}}
               \raise2pt\hbox{\rlap{\hskip-9.5pt \vrule height.4pt depth0pt
                width6.2pt}}
               \lower6pt\hbox{\rlap{\vrule height4.5pt}}}}
\def\je{\smash{\raise2pt\hbox{\co \rlap{\rlap{\char'005}
                \phantom{\char'007}}}\raise6pt\hbox{\rlap{\vrule height5pt}}
               \raise2pt\hbox{\rlap{\hskip-9.5pt \vrule height.4pt depth0pt
                width6.2pt}}}}
\def\er{\smash{\raise2pt\hbox{\co \rlap{\rlap{\phantom{\char'005}} \char'007}}
               \raise2pt\hbox{\rlap{\hskip4pt \vrule height0.4pt depth0pt
                width5.7pt}}
               \lower6pt\hbox{\rlap{\vrule height4.5pt}}}}
\def\+{\smash{\lower6pt\hbox{\rlap{\vrule height17pt}}
                \raise2pt\hbox{\rlap{\hskip-9pt \vrule height.4pt depth0pt
                width18.7pt}}}}
\def\hor{\smash{\raise2pt\hbox{\rlap{\hskip-9.5pt \vrule height.4pt depth0pt
                width19.2pt}}}}
\def\ver{\smash{\lower6pt\hbox{\rlap{\vrule height17pt}}}}

\def\textcross{\ \smash{\lower4pt\hbox{\rlap{\hskip4.15pt\vrule height14pt}}
                \raise2.8pt\hbox{\rlap{\hskip-3pt \vrule height.4pt depth0pt
                width14.7pt}}}\hskip12.7pt}
\def\textelbow{\ \hskip.1pt\smash{\raise2.8pt%
                \hbox{\co \hskip 4.15pt\rlap{\rlap{\char'005} \char'007}
                \lower6.8pt\rlap{\vrule height3.5pt}
                \raise3.6pt\rlap{\vrule height3.5pt}}
                \raise2.8pt\hbox{%
                  \rlap{\hskip-7.15pt \vrule height.4pt depth0pt width3.5pt}%
                  \rlap{\hskip4.05pt \vrule height.4pt depth0pt width3.5pt}}}
                \hskip8.7pt}

\newcommand{\cellsize}{22}
\newlength{\cellsz} 
\newsavebox{\cell}
\sbox{\cell}{\begin{picture}(\cellsize,\cellsize)
\put(0,0){\line(1,0){\cellsize}}
\put(0,0){\line(0,1){\cellsize}}
\put(\cellsize,0){\line(0,1){\cellsize}}
\put(0,\cellsize){\line(1,0){\cellsize}}
\end{picture}}
\newcommand\cellify[1]{\def\thearg{#1}\def\nothing{}%
\ifx\thearg\nothing
\vrule width0pt height\cellsz depth0pt\else
\hbox to 0pt{\usebox{\cell} \hss}\fi%
\vbox to \cellsz{
\vss
\hbox to \cellsz{\hss$#1$\hss}
\vss}}
\newcommand\tableau[1]{\vtop{\let\\\cr
\baselineskip -16000pt \lineskiplimit 16000pt \lineskip 0pt
\ialign{&\cellify{##}\cr#1\crcr}}}
\begin{document}
\pagestyle{plain}

\title{When is a Schubert Variety Gorenstein?}
\subjclass[2000]{14M15; 14M05, 05E99}
\author{Alexander Woo}
\address{Department of Mathematics, Kerr Hall, One Shields Ave., University of California, Davis, CA 95616, USA}
\email{awoo@math.ucdavis.edu}

\author{Alexander Yong}
\address{Department of Mathematics, University of Minnesota, Minneapolis, MN 55455, USA,\indent{\itshape and}
        The Fields Institute, 222 College Street, Toronto, Ontario,  
        M5T 3J1, Canada}

\email{ayong@math.umn.edu, ayong@fields.utoronto.ca}

\date{November 20, 2005}

\maketitle

\section{Introduction}
The main goal of this paper is to give an explicit combinatorial
characterization of which
Schubert varieties in the complete flag variety are Gorenstein. 

Let $\mathrm{Flags}({\mathbb C}^{n})$ denote the variety of complete
flags $F_{\bullet}:\langle 0\rangle \subseteq
F_1\subseteq\ldots\subseteq F_n={\mathbb C}^n$. Fix a basis
$e_1,e_2,\ldots,e_n$ of $\mathbb{C}^n$ and let $E_{\bullet}$ be the
{\em anti}-canonical reference flag $E_{\bullet}$, that is, the flag
where $E_{i}=\langle e_{n-i+1},e_{n-i+2},\ldots, e_n\rangle$.  For
every permutation $w$ in the symmetric group $S_n$, there is the {\bf
Schubert variety}
\[X_{w}=\big\{F_{\bullet}\mid {\rm dim}(E_i\cap F_j)\geq \#\{k\geq n-i+1, w(k)\leq j\}\big\}.\]
These conventions have been arranged so that the codimension of
$X_{w}$ is $\ell(w)$, that is, the length of any expression for $w$ as a
product of simple reflections $s_i = (i\leftrightarrow i+1)$.  

The Gorenstein property gives a well-known measurement of how far an
algebraic variety is from being smooth; all smooth varieties are
Gorenstein, while all Gorenstein varieties are Cohen-Macaulay.  In
general, a variety is {\bf Gorenstein} if it is Cohen-Macaulay and its
canonical sheaf is a line bundle.  (Throughout this paper we freely
identify vector bundles and their sheaves of sections for
convenience.)  Recall that on a {\em smooth} variety $X$, the {\bf
canonical sheaf}, denoted $\omega_\X$, is $\bigwedge^{{\rm
dim}(\X)}\Omega_{\X}$, where $\Omega_{\X}$ is the cotangent bundle of
$\X$. For a possibly singular but normal variety $\X$, the canonical sheaf
is the pushforward of the canonical sheaf $\omega_{\X_{\rm smooth}}$
on the smooth part $X_{\rm smooth}$ of $\X$ under the inclusion map.
Since every Schubert variety is normal~\cite{DeCon-Lak,RR} and
Cohen-Macaulay~\cite{Ram85}, the remarks above suffice to define
Gorensteinness in the context of this paper. In Section~2.1,
we will give the more commonly seen local definition of Gorensteinness.
However, combining the above definition together with the results of
Ramanathan~\cite{Ram85,Ram87} is what provides our starting
point for determining which Schubert varieties are Gorenstein.

Smoothness and Cohen-Macaulayness of Schubert varieties have been
extensively studied in the literature; see, for
example,~\cite{BL00,Ram85} and the references therein. While all
Schubert varieties are Cohen-Macaulay, very few Schubert
varieties are smooth.  (See the table at the end of this
Introduction.)  Explicitly, $X_w$ is smooth if and only if $w$ is
``$1324$-pattern avoiding'' and ``$2143$-pattern avoiding''
\cite{LS90}; we give more details on pattern avoidance below.

Our main result (Theorem~\ref{thm:main}) gives an explicit
combinatorial characterization of which Schubert varieties are
Gorenstein similar to the above smoothness criteria.  This answers a
question raised by M.~Brion and S.~Kumar and passed along to us by
A.~Knutson; see also~\cite[p.~88]{Ram87}.  Our answer uses a
generalized notion of pattern avoidance that we introduce.

To describe our ideas in a simpler case, we first compare the
classical smoothness criterion with a characterization of which
Schubert varieties in the Grassmannian $\mathrm{Gr}(\ell,n)$ of
$\ell$-planes in ${\mathbb C}^n$ are Gorenstein.  (This is a special
case of our main result, as we will explain in
Section~\ref{subsection:partial}.)  Schubert varieties $\X_{\lambda}$
of $\mathrm{Gr}(\ell,n)$ are indexed by partitions $\lambda$ sitting
inside an $\ell\times (n-\ell)$ rectangle.\footnote{Consistent with
our convention on Schubert varieties in
$\mathrm{Flags}(\mathbb{C}^n)$, we index these Schubert varieties so
that $\left|\lambda\right|$ is the codimension of $X_\lambda$.} The
smooth Schubert varieties are those indexed by partitions $\lambda$
whose complement in $\ell\times (n-\ell)$ is a rectangle, as explained
in, for example,~\cite{BL00} and the references therein.  For example,
$\lambda=(7,7,2,2,2)$ indexes a smooth Schubert variety in ${\rm
Gr}(5,12)$.
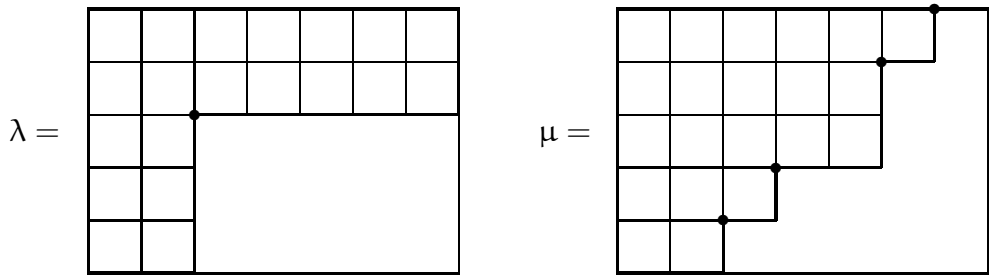
\[
\begin{picture}(340,100)
\put(0,0){\framebox(140,100)}
\thicklines
\put(0,0){\line(1,0){40}}
\put(40,0){\line(0,1){60}}
\put(40,60){\line(1,0){100}}
\put(140,60){\line(0,1){40}}
\put(-30,50){$\lambda=$}
\put(40,60){\circle*{4}}
\put(0,0){\line(0,1){100}}
\put(0,100){\line(1,0){140}}
\thinlines
\put(0,20){\line(1,0){40}}
\put(0,40){\line(1,0){40}}
\put(0,60){\line(1,0){40}}
\put(0,80){\line(1,0){140}}
\put(20,0){\line(0,1){100}}
\put(40,60){\line(0,1){40}}
\put(60,60){\line(0,1){40}}
\put(80,60){\line(0,1){40}}
\put(100,60){\line(0,1){40}}
\put(120,60){\line(0,1){40}}

\put(200,0){\framebox(140,100)}
\thicklines
\put(200,0){\line(1,0){40}}
\put(240,0){\line(0,1){20}}
\put(240,20){\line(1,0){20}}
\put(260,20){\line(0,1){20}}
\put(260,40){\line(1,0){40}}
\put(300,40){\line(0,1){40}}
\put(300,80){\line(1,0){20}}
\put(320,80){\line(0,1){20}}
\put(320,100){\circle*{4}}
\put(170,50){$\mu=$}
\put(240,20){\circle*{4}}
\put(260,40){\circle*{4}}
\put(300,80){\circle*{4}}
\put(200,0){\line(0,1){100}}
\put(200,100){\line(1,0){140}}
\thinlines
\put(220,0){\line(0,1){100}}
\put(240,20){\line(0,1){80}}
\put(260,40){\line(0,1){60}}
\put(280,40){\line(0,1){60}}
\put(300,80){\line(0,1){20}}
\put(200,20){\line(1,0){40}}
\put(200,40){\line(1,0){60}}
\put(200,60){\line(1,0){100}}
\put(200,80){\line(1,0){100}}
\end{picture}
\]

Alternatively, smooth Schubert varieties are those with at most one
inner corner.  View the lower border of partition as a lattice path
from the lower left-hand corner to the upper right-hand corner of
$\ell\times (n-\ell)$; an {\bf inner corner} is then a lattice point
on this path with lattice points of the path both directly
below and directly to the right of it. The inner corners for the
partitions $\lambda$ and $\mu$ above are marked by ``dots''.

Therefore, the partition $\mu=(6,5,5,3,2)$ above does not index a
smooth Schubert variety.  However, it does index a Gorenstein Schubert
variety; in general, a Grassmannian Schubert variety $\X_{\mu}$ is
Gorenstein if and only if all of the inner corners of $\mu$ lie on the
same antidiagonal. We mention that this condition can also be derived
from \cite[(5.5.5)]{svanes}.

In order to state our main result for $\mathrm{Flags}({\mathbb
C}^{n})$, we will need some preliminary definitions.  First we
associate a Grassmannian permutation to each descent of a permutation
$w$.  Let $d$ be a {\bf descent} of $w$, which is an index such that
$w(d)>w(d+1)$.  Now write $w$ in one-line notation as $w(1)w(2)\cdots
w(n)$, and construct a subword $v_{d}(w)$ of $w$ by concatenating the
right-to-left minima of the segment strictly to the left of $d+1$ with
the left-to-right maxima of the segment strictly to the right of
$d$. In particular, $v_{d}(w)$ will necessarily include $w(d)$ and
$w(d+1)$.  Let ${\widetilde v}_{d}(w)$ denote the {\bf flattening} of
$v_{d}(w)$, which is defined to be the unique permutation whose
entries are in the same relative position as those of $v_{d}(w)$.

\begin{Example} Let $w=314972658\in S_{9}$.
This permutation has descents at positions 1, 4, 5 and 7.  We see
that $v_{1}(w)=3149$, $v_{4}(w)=14978$, $v_{5}(w)=147268$, and
$v_{7}(w)=12658$, so therefore ${\widetilde v}_{1}(w)=2134$,
${\widetilde v}_{4}(w)=12534$, ${\widetilde v}_{5}(w)=135246$, and
${\widetilde v}_{7}(w)=12435$.
\end{Example}

By construction, ${\widetilde v}_{d}(w)\in S_m$ is a {\bf Grassmannian
permutation}, meaning that it has a unique descent at some position we
denote $e$.  For any Grassmannian permutation $w\in S_m$ with its
unique descent at $e$, let $\lambda(w)\subseteq e\times (m-e)$ denote
the associated partition.  The partition $\lambda(w)$ is the one whose
lower border is obtained by drawing a lattice path which starts at the
lower left corner of $e\times(m-e)$ and continues by a unit horizontal
line segment at step $i$ (for some $i\in\{1,\ldots,m\}$) if $i$
appears strictly after position $e$ (or, in other words, if
$w^{-1}(i)>e$), and a unit vertical line segment otherwise.  For
example, the Grassmannian permutation~$w= 3589\ 11 \ \mid\ 12467\ 10\
12$ corresponds to the partition $\lambda(w)=\mu=(6,5,5,3,2)$ depicted
above.  Now, given an inner corner of a partition $\lambda(w)$, let
its {\bf inner corner distance} be the sum of the distances from the
top and left edges of the rectangle $e\times (m-e)$ to the inner
corner. For example, in $\mu$ above, all the inner corner distances
equal~6.  Furthermore, suppose that $\lambda(w)$ has all its inner
corners on the same antidiagonal; this is equivalent to requiring that
the inner corner distance be the same for all inner corners.  In this
case we call this common inner corner distance ${\mathfrak I}(w)$; if
there are no inner corners, we set $\mathfrak{I}(w)=0$ by convention.
For our example permutation~$w$, $\mathfrak{I}(w)=6$.

Next we proceed to define {\bf Bruhat-restricted pattern avoidance}.
Recall that, classically, for $v\in S_\ell$ and $w\in S_n$, with
$\ell\leq n$, an {\bf embedding of $v$ into $w$} is a sequence of
indices $i_1<i_2<\cdots<i_{\ell}$ such that, for all $1\leq
a<b\leq\ell$, $w(i_a)>w(i_b)$ if and only if $v(a)>v(b)$.  Then the
classical definition of pattern avoidance is that {\bf $w$ pattern
avoids $v$} if there are no embeddings of $v$ into $w$.

Now recall the {\bf Bruhat order} $\succ$ on $S_n$.  First we say that
$w(i\leftrightarrow j)$ {\bf covers} $w$ if $i<j$, $w(i)<w(j)$, and,
for each $k$ with $i<k<j$, either $w(k)<w(i)$ or $w(k)>w(j)$; then the
Bruhat order is the transitive closure of this covering relation.  The
Bruhat order is graded by the length of a permutation, and one can
check that $v$ can cover $w$ only if $\ell(v)=\ell(w)+1$.

Given a permutation $v\in S_\ell$, let ${\mathcal T}_{v}=\{(m_1
\leftrightarrow n_1), \ldots, (m_k \leftrightarrow n_k)\}$ be a set of
{\bf Bruhat transpositions} in $v$, by which we mean a subset of
transpositions such that $v\cdot(m_j \leftrightarrow n_j)$ covers $v$
in the Bruhat order.  We define a {\bf ${\mathcal T}_{v}$-restricted
embedding} of $v$ into $w$ to be an embedding of $v$ into $w$ such
that $w\cdot(i_{m_j} \leftrightarrow i_{n_j})$ covers $w$ for all
$(m_j \leftrightarrow n_j)\in {\mathcal T}_{v}$. Then we say that {\bf
$w$ pattern avoids $v$ with Bruhat restrictions ${\mathcal T}_{v}$} if
there are no ${\mathcal T}_{v}$-restricted embeddings of $v$ into $w$.

Now we are ready to state our combinatorial characterization of which
Schubert varieties in $\mathrm{Flags}({\mathbb C}^{n})$ are
Gorenstein:

\begin{Theorem}
\label{thm:main}
Let $w\in S_n$. The Schubert variety $X_w$ is Gorenstein if and only if
\begin{itemize}
\item for each descent $d$ of $w$, $\lambda({\widetilde v}_{d}(w))$ has all
of its inner corners on the same antidiagonal, and
\item the permutation $w$ pattern avoids both $31524$ and $24153$ with
Bruhat restrictions \linebreak $\{(1\leftrightarrow 5),(2\leftrightarrow 3)\}$
and $\{(1\leftrightarrow 5),(3\leftrightarrow 4)\}$ respectively.
\end{itemize}
\end{Theorem}

In comparing the smoothness characterization of \cite{LS90} with
Theorem~\ref{thm:main}, considering our description of the
Grassmannian case allows one to check that the $1324$-pattern avoidance
condition of the former implies the ``inner corner condition'' of the
latter. It is also easy to see that the $2143$-pattern avoidance
condition of the former implies both of the Bruhat-restricted pattern
avoidance conditions of the latter. We mention that
Fulton~\cite{Fulton:Duke92} has characterized $2143$-pattern avoidance
in terms of the essential set of a permutation.  A similar
characterization can be given for the Bruhat-restricted pattern
avoidance conditions of Theorem~\ref{thm:main}.

\medskip
\begin{Example}
The permutation $w=\underline{3} 7 \underline{1} 4 \underline{8}
\underline{2}6\underline{5}\in S_8$ has descents at positions $2$, $5$ and $7$ and we have 
\[{\widetilde v}_{2}(w)=24135, {\widetilde v}_{5}(w)=13524, \mbox{ and } {\widetilde v}_{7}(w)=1243.\]
Hence one checks that $w$ satisfies the inner corner condition with
\[{\mathfrak I}({\widetilde v}_{2}(w))=2, \ {\mathfrak I}({\widetilde v}_{5}(w))=2,
\mbox{ and } {\mathfrak I}({\widetilde v}_{7}(w))=1.\] The Schubert
variety $X_{w}$ is Gorenstein, since there are no forbidden $31524$
and $24153$ patterns with Bruhat restrictions $\{(1\leftrightarrow
5),(2\leftrightarrow 3)\}$ or $\{(1\leftrightarrow 5),
(3\leftrightarrow 4)\}$ respectively. Note that the underlined subword
of $w$ is a $31524$-pattern, but since $w(1\leftrightarrow 8)$ does
not cover $w$, it does not prevent $X_w$ from being Gorenstein.
\end{Example}

By combining Theorem~\ref{thm:main} with the descriptions of the
singularities along the ``maximal singular locus'' of a Schubert
variety $X_w$ given in~\cite{Cortez, manivel2}, we obtain the following
geometric corollary.

\begin{Corollary}
\label{cor:maxsingloc}
A Schubert variety $X_w$ is Gorenstein if and only if it is Gorenstein
along its maximal singular locus.
\end{Corollary}

In other words, Corollary~\ref{cor:maxsingloc} states that a Schubert
variety is Gorenstein if and only if its ``smoothest'' singularities
(those at the generic points of the irreducible components of the
singular locus) are Gorenstein.

	We now describe the canonical sheaf of a Gorenstein Schubert
variety in terms of the Borel-Weil construction of line bundles.  Let
$T\cong(\mathbb{C}^*)^{n-1}$ be the subgroup of invertible diagonal
matrices of determinant $1$ in $\mathrm{SL}_n(\mathbb{C})$; the
Borel-Weil construction associates to each integral weight
$\alpha\in\mathrm{Hom}(T,\mathbb{C}^*)$ a line bundle
$\mathcal{L}_\alpha$.  Let ${\mathcal L}_{\alpha}\big|_{X_{w}}$ denote
the restriction of this line bundle to $X_w$.  We will write weights
additively in terms of the $\mathbb{Z}$-basis of fundamental weights
$\Lambda_r$, defined by $\Lambda_r\left(\left[\begin{array}{ccc} t_1 &
& 0 \\ & \ddots & \\ 0 & & t_n\end{array}\right]\right)=t_1\cdots
t_r$.

\begin{Theorem}
\label{thm:second_main}
If $X_{w}$ is Gorenstein, then 
$\omega_{X_{w}}\cong {\mathcal L}_{\alpha}\big|_{X_w}$
where $\alpha=\sum_{r=1}^{n-1} \widetilde{\alpha}_r\Lambda_{n-r}$ and
\begin{equation}
\label{eqn:the_soln'}
\widetilde{\alpha}_r = \left\{
\begin{array}{cc}
-2+{\mathfrak I}({\widetilde v}_{r}(w)) & \mbox{ if $r$ is a descent} \\ 
-2 & \mbox{ otherwise.}
\end{array}
\right.
\end{equation}
\end{Theorem}

	The proofs of Theorems~\ref{thm:main}
and~\ref{thm:second_main} as well as Corollary~\ref{cor:maxsingloc}
will be given in Section~2. In Section~3, we end with a number of
remarks and applications.

	Further study of the relationship between the 
geometry of Gorensteinness of Schubert varieties and related combinatorics 
should have potential. We conclude this introduction with some open 
problems and suggestions for further work. 
The most natural is:
\begin{Problem}
\label{problem:1}
Give analogues of Theorems~\ref{thm:main} and~\ref{thm:second_main} for
generalized flag varieties corresponding to Lie groups other than ${\rm GL}_{n}({\mathbb C})$.
\end{Problem}
\noindent
We expect that the methods given in this paper will extend to solve Problem~\ref{problem:1}. It is not difficult to use Theorem~\ref{thm:main} to derive an analogue of Theorem~\ref{thm:main} for the case of the odd 
orthogonal groups
${\rm SO}(2n+1,{\mathbb C})$. However, we have found the 
combinatorial analysis required to be more intricate in general. 
Consequently, in the interest of brevity, 
we plan to discuss our investigations for the other Lie types in a subsequent paper.

	It should also be interesting to determine the ``maximal
non-Gorenstein locus'' of a non-Gorenstein Schubert variety: Let $X$
be a variety that is Cohen-Macaulay but not Gorenstein; since the rank
of any coherent sheaf on $X$ is upper semicontinuous (see, for
example,~\cite[III.12.7.2]{Hartshorne}), the canonical sheaf has rank
strictly greater than 1 at some non-trivial closed subvariety.  This
subvariety then consists of all points of $X$ at which $X$ is not
Gorenstein by the local definition.  Since the canonical sheaf of a
Schubert variety is $B_{-}$-equivariant for the subgroup
$B_{-}\subseteq {\rm GL}_{n}({\mathbb C})$ of lower triangular
matrices, this subvariety is a union of Schubert varieties contained
in $X_w$. Therefore we ask:
\begin{Problem}
Give a combinatorial characterization for
the minimal $v$ in the Bruhat order for which
$X_{w}$ is not Gorenstein at $X_v$.
\end{Problem}

In view of Corollary~\ref{cor:maxsingloc}, it is natural
to propose the following answer:
\begin{Conjecture}
\label{conj:nongorloc}
The maximal non-Gorenstein locus of $X_w$ is the union of those 
Schubert varieties $X_v$ in the maximal singular locus of $X_w$
for which the generic point is not Gorenstein in $X_w$.
\end{Conjecture}

One can give a combinatorial rule characterizing the set of $X_v$
appearing in Conjecture~\ref{conj:nongorloc} using the explicit
description of the singular locus of Schubert
varieties~\cite{billey.warrington, Cortez, Gasharov, KLR, LS90,
manivel1, manivel2} and facts mentioned in the proof of
Corollary~\ref{cor:maxsingloc}.

A geometric explanation was recently given in \cite{BilBrad} for the
appearance of pattern avoidance in characterizations of smooth
Schubert varieties.  However, this explanation does not have an
obvious modification to take into account Bruhat-restrictions.  This
leads to the following:
\begin{Problem}
Give a geometric explanation of Bruhat-restricted pattern avoidance
which explains its appearance in Theorem~\ref{thm:main}.
\end{Problem}

	Lastly, for those interested in combinatorial enumeration:
\begin{Problem}
Give a combinatorial formula (for example, a generating series)
computing the number of Gorenstein Schubert varieties in ${\rm
Flags}({\mathbb C}^n)$.
\end{Problem}
\noindent
Using the methods of this paper, we computed the number of
Gorenstein Schubert varieties in ${\rm Flags}({\mathbb C}^{n})$ for some
small values of $n$ (see below). 
We compare this to the number of smooth Schubert 
varieties computed using the result of~\cite{LS90} 
(by the recursive formulas found in \cite{Bona,Stankova}).

\begin{table}[h]
\begin{center}
\begin{tabular}{|l|l|l|l|}
\hline
$n$ & $n! = \#\mbox{ Cohen-Macaulay $X_w$}$ & $\#\mbox{ Gorenstein $X_w$}$  & $\#\mbox{ Smooth $X_w$}$ \\ \hline\hline

1 & 1 & 1 & 1 \\ \hline
2 & 2 & 2 & 2 \\ \hline
3 & 6 & 6 & 6 \\ \hline
4 & 24 & 24 & 22  \\ \hline
5 & 120 & 116 & 88 \\ \hline
6 & 720 & 636 & 366 \\ \hline
7 & 5040 & 3807 & 1552 \\ \hline
8 & 40320 & 24314 & 6652 \\ \hline
9 & 362880 & 163311 & 28696 \\ \hline
\end{tabular}
\end{center}
\end{table}

We are very grateful to M.~Brion, A.~Knutson and S.~Kumar for bringing
the problem addressed by Theorem~\ref{thm:main} to our attention, for
outlining the argument used in Section~2.1, and for many other
suggestions.  We also thank A.~Bertram, S.~Billey, A.~Buch, A.~Cortez,
R.~Donagi, S.~Fomin, M.~Haiman, R.~MacPherson, E.~Miller, R.~Stanley,
B.~Sturmfels, J.~Tymoczko, and an anonymous referee for discussion and
remarks on earlier drafts.  This work was partially completed while
the two authors were in residence at the Park City Mathematics
Institute program on ``Geometric Combinatorics'' during July 2004.

\section{Proof of Theorems~\ref{thm:main} and~\ref{thm:second_main}}

\subsection{Geometry to combinatorics}

First we explain the algebraic definition of Gorensteinness and reduce
the algebro-geometric problem of determining when a Schubert variety
is Gorenstein to a problem in linear algebra; we will then solve this
linear algebra problem combinatorially.  This reduction to linear
algebra appears to be folklore (and was told to us by M.~Brion,
A.~Knutson and S.~Kumar); we could not locate an explicit reference
for it in the literature.  Therefore, we include an argument for the
sake of completeness.  While we treat only
$\mathrm{Flags}(\mathbb{C}^n)$ explicitly, the arguments of this
section generalize easily to all semi-simple Lie groups with the
substitution of the appropriate Monk-Chevalley
formula~\cite{Chevalley}.  We found~\cite{brion:book} an excellent
resource for facts about the geometry of Schubert varieties.

A local ring $(R,\mathfrak{m}, \Bbbk)$ is said to be {\bf Cohen-Macaulay}
if $\mathrm{Ext}_R^i(\Bbbk, R)=0$ for $i\leq\mathrm{dim} R$; it is {\bf
Gorenstein} if, in addition, $\mathrm{dim}_{\Bbbk}
\mathrm{Ext}_R^{\mathrm{dim} R}(\Bbbk, R)=1$.  A variety is Cohen-Macaulay
(respectively Gorenstein) if the local ring at every point is
Cohen-Macaulay (respectively Gorenstein).  Using the Kozsul complex on
a regular sequence, one can show that every regular local ring is
Gorenstein; hence smooth varieties are Gorenstein.  See
\cite{Bruns-Herzog} for details.

One might naively expect that, in order to check if a Schubert variety
is Gorenstein, one would need to check if it is Gorenstein at all, or
at least some, of its points.  However, the alternative equivalent
definition of the Gorenstein property alluded to in the introduction,
which is based on Grothendieck duality theory (see \cite{HartsRes} or
\cite{Altman-Kleiman}), allows for a different approach using the
global geometry of Schubert varieties.  Each projective variety has a
dualizing complex (of sheaves) which plays a role analogous to that of
the canonical bundle $\omega_X$, defined as the top exterior power of
the cotangent bundle $\bigwedge^{{\rm dim}(X)}\Omega_{X}$, of a smooth
variety in Serre duality.  A connected projective variety is
Cohen-Macaulay if and only if the dualizing complex is a sheaf, and
Gorenstein if and only if the dualizing sheaf is locally free of
rank one.  For a normal, Cohen-Macaulay variety, one can realize the
dualizing sheaf as the pushforward of the canonical sheaf
$\omega_{X_\mathrm{smooth}}$ of the smooth part $X_{\mathrm{smooth}}$
under the inclusion map.  As mentioned in the introduction, all
Schubert varieties are known to be normal~\cite{DeCon-Lak,RR} and
Cohen-Macaulay~\cite{Ram85}, so we can then use the calculation of the
canonical sheaf of Schubert varieties by Ramanathan~\cite{Ram85,Ram87}
to determine which Schubert varieties are Gorenstein.

We now need some standard definitions which can be found
in~\cite[II.6]{Hartshorne}.  Let $\mathrm{Cl}(X_w)$ denote the Weil
divisor class group of $X_w$; its elements are linear equivalence
classes $\left[ Z\right]$ of formal sums of codimension~1 subvarieties
$Z$ of $X_w$.  There is a natural group homomorphism $\mathrm{div}:
\mathrm{Pic}(X_w) \rightarrow \mathrm{Cl}(X_w)$, where
$\mathrm{Pic}(X_w)$ is the group of isomorphism classes of line
bundles under tensor product.  On a Schubert variety $X_w$ (or, in
general, any normal irreducible variety over a field), $\mathrm{div}$
is injective and its image in $\mathrm{Cl}(X_w)$ is the {\bf Cartier
class group} $\mathrm{CaCl}(X_{w})$.  (This is an unorthodox
definition of the Cartier class group, but for convenience we have
identified it with its isomorphic image in the Weil class group.)  For
smooth varieties, $\mathrm{div}$ is an isomorphism, so
$\mathrm{CaCl}=\mathrm{Cl}$.

We now proceed to describe explicitly $\mathrm{Cl}(X_w)$ and
$\mathrm{CaCl}(X_w)$.  The Schubert variety $X_w$ is the disjoint
union of the open Schubert cell $X^\circ_w$ (which is isomorphic to
the affine space ${\mathbb C}^{\binom{n}{2}-\ell(w)}$) together with
the codimension~1 subvarieties $X_v$ for $v$ covering $w$ in the
Bruhat order.  Therefore, by repeatedly applying
\cite[Prop. II.6.5]{Hartshorne}, we see that $\mathrm{Cl}(X_w)$ is
freely generated (as an abelian group) by $\left[ X_v \right]$ for $v$
covering $w$.

\smallskip
To describe $\mathrm{CaCl}(X_w)$, we will need the Chow group
$A_*(\mathrm{Flags}(\mathbb{C}^{n}))$ of the flag variety, whose
elements are rational equivalence classes $\left[Z\right]$ of
subvarieties $Z$ of $\mathrm{Flags}(\mathbb{C}^{n})$; see for
example~\cite[Ch.~1]{FultonIT}.  Since
$\mathrm{Flags}(\mathbb{C}^{n})$ is smooth, the Chow ring
$A^*(\mathrm{Flags}(\mathbb{C}^{n}))$ is by definition equal as
abelian groups to
$A_*(\mathrm{Flags}(\mathbb{C}^{n}))$~\cite[8.3]{FultonIT}.  The
graded pieces
$A_d(\mathrm{Flags}(\mathbb{C}^n))=A^{\binom{n}{2}-d}(\mathrm{Flags}(\mathbb{C}^{n}))$
are freely generated by the classes $\left[X_v\right]$ of the Schubert
varieties $X_v$ of dimension $d$, which are precisely those for which 
$d=\binom{n}{2}-l(v)$.  Therefore, the natural map $\iota_*:
\mathrm{Cl}(X_w) \rightarrow
A_{\binom{n}{2}-\ell(w)-1}(\mathrm{Flags}(\mathbb{C}^n))$ induced by
the inclusion \linebreak $\iota: X_w \rightarrow \mathrm{Flags}(\mathbb{C}^n)$ is
injective.  Note that, by definition,
$A_{\binom{n}{2}-1}(\mathrm{Flags}(\mathbb{C}^{n}))=\mathrm{Cl}(\mathrm{Flags}(\mathbb{C}^{n}))$.

\smallskip
It is known~\cite[Prop. 6]{Mathieu} that every line bundle on a
Schubert variety is the restriction of a line bundle on
$\mathrm{Flags}(\mathbb{C}^{n})$.  Furthermore, for a line bundle
$\mathcal{L}$ on $\mathrm{Flags}(\mathbb{C}^{n})$, general facts of
intersection theory~\cite[Ch. 2]{FultonIT} tell us that
$\iota_*(\mathrm{div}(\mathcal{L}\big|_{X_w}))=\mathrm{div}(\mathcal{L})\cdot\left[X_w\right]$,
where the right hand side is a product in
$A^*(\mathrm{Flags}(\mathbb{C}^{n}))$.  Therefore, since
$\mathrm{CaCl}(\mathrm{Flags}(\mathbb{C}^{n}))=\mathrm{Cl}(\mathrm{Flags}(\mathbb{C}^{n}))$
is generated by $\left\{\left[X_{(r\leftrightarrow
r+1)}\right]\right\}_{r=1}^{n-1}$,
$\iota_*(\mathrm{CaCl}(X_w))\subseteq {\rm A}^{*}({\rm Flags}({\mathbb
C}^n))$ is generated by $\left\{\left[X_{(r\leftrightarrow
r+1)}\right]\cdot\big[X_w\big]\right\}_{r=1}^{n-1}$.

\smallskip
By Monk's formula~\cite{Monk}, 
$$\left[X_{(r\leftrightarrow r+1)}\right]\cdot\big[X_w\big]
=\sum_{\plainstack{a\leq r<b}{\ell(w(a\leftrightarrow b))=\ell(w)+1}}
\left[ X_{w(a\leftrightarrow b)}\right],$$ so $\mathrm{CaCl}(X_w)$ is
generated by these classes for $1\leq r\leq n-1$.  (We can drop the
$\iota_*$ since it is an injection.)

\smallskip
Since Schubert varieties are Cohen-Macaulay~\cite{Mus-Ses,DeCon-Lak,Ram85},
a Schubert variety $X_w$ is Gorenstein if and only if its canonical
sheaf $\omega_{X_w}$ is a line bundle.  By results of
Ramanathan~\cite[Thm. 4.2]{Ram87},
the canonical sheaf of $X_w$ is
$$\omega_{X_w}=\mathcal{L}_{-\rho}\mid_{X_w} \otimes \
\mathcal{I}(\partial X_w),$$ where $\mathcal{L}_{-\rho}\mid_{X_w}$ is
the restriction to $X_w$ of the line bundle associated to the weight
$-\rho=-\sum_{r=1}^{n-1}\Lambda_{r}$ by the Borel-Weil construction,
and $\mathcal{I}(\partial X_w)$ is the ideal sheaf of the complement
of $X^\circ_w$, or equivalently, the ideal sheaf of the reduced
subscheme $\bigcup_v X_v$ where $v$ ranges over all permutations
covering $w$ in the Bruhat order.  Since
$\mathcal{L}_{-\rho}\big|_{X_w}$ is a line bundle and $\mathrm{Pic}$
is a group, $\omega_{X_w}$ is a line bundle if and only if
$\mathcal{I}(\partial X_w)$ is a line bundle.  However, the ideal
sheaf of a reduced codimension~1 subscheme $Y$ is a line bundle if and
only if $\left[Y\right]$ is a Cartier divisor, in which case
$\mathrm{div}(\mathcal{I}(Y))=-\left[Y\right]$; see, for example,
\cite[II.6]{Hartshorne}.  Therefore, $X_w$ is Gorenstein if and only
if $$\left[\partial X_w\right]=\sum_{\plainstack{v\succ
w}{\ell(v)=\ell(w)+1}} \left[\X_v\right]\in \mathrm{CaCl}(X_{w}).$$

\pagebreak
Hence, by our previous calculation of $\mathrm{CaCl}(X_w)$ as a
subgroup of $\mathrm{Cl}(X_w)$, we obtain the following:

\begin{Proposition}
The Schubert variety $X_w$ is Gorenstein
if and only if there exists
an integral solution $(\alpha_1,\ldots,\alpha_{n-1})$ to
\begin{equation}
\label{eqn:final_system}
\sum_{r=1}^{n-1}\alpha_r \left(\sum_{\plainstack{a\leq r<b}{\ell(w(a\leftrightarrow b))=\ell(w)+1}}\left[ X_{w(a\leftrightarrow b)}\right] \right)= \sum_{\plainstack{v=w(a\leftrightarrow b)}{\ell(v)=\ell(w)+1}} \left[ X_v\right]\in {\rm CaCl}(X_w).
\end{equation}
\end{Proposition}

\smallskip
As an aside, a variety is said to be {\bf locally factorial} if the
local ring at every point is a unique factorization domain.  It is
well known (see \cite[Prop. II.6.11]{Hartshorne} or
\cite[2.1]{FultonIT}) that a normal variety is factorial if and only
if $\mathrm{div}$ is an isomorphism.  Therefore, factorial Schubert
varieties can be characterized using the following proposition.

\begin{Proposition}
The Schubert variety $X_w$ is factorial if and only if the classes
$$\left\{\sum_{\plainstack{a\leq r<b}{\ell(w(a\leftrightarrow
b))=\ell(w)+1}}\left[ X_{w(a\leftrightarrow
b)}\right]\right\}_{r=1}^{n-1}$$ span the free abelian group generated
by $$\left\{\left[X_v\right] \mid v=w(a\leftrightarrow b),
\ell(v)=\ell(w)+1\right\}.$$
\end{Proposition}

\smallskip
Recently, M. Bosquet-M\'elou and S. Butler~\cite{BM-Butler} have used
this proposition to give a characterization of locally factorial
Schubert varieties in terms of Bruhat-resticted pattern avoidance.
This solves a conjecture that we had distributed during the
preparation of this article.
\vspace{-.2in}
\subsection{Interlude: a diagrammatic formulation and two sample problems}

	Although it is not used in our proof below, let us 
give a diagrammatic formulation of the above linear algebra problem 
(\ref{eqn:final_system}) that the reader may find useful. 
 
\smallskip
	Label $n$ columns by the values $w(1), w(2), \ldots, w(n)$ of
a permutation $w\in S_n$. Draw horizontal bars between the midpoints
of columns $i$ and $j$ if and only if $w(i\leftrightarrow j)$ covers
$w$ in the Bruhat order. Now draw vertical bars between columns $i$ and $i+1$ for $1\leq
i\leq n-1$. Then a solution to (\ref{eqn:final_system}) is equivalent
to an assignment $(\alpha_1,\ldots,\alpha_{n-1})\in {\mathbb Z}^{n-1}$
of integers to the vertical bars (from left to right respectively)
such that, for each horizontal bar, the sum of the assignments to the
vertical bars that it crosses equals~1.

	We encourage the reader to try out the following two sample problems; 
answers are at the bottom of the page\footnote{The problem on the
left is solved by $(\alpha_1,\alpha_2,\alpha_3,\alpha_4,
\alpha_5,\alpha_6)=(-1,0,1,1,-1,1)$ while the problem
on the right has no solution.}:
\[
\begin{picture}(320,200)
\put(0,190){$6$}
\put(20,190){$3$}
\put(40,190){$1$}
\put(60,190){$4$}
\put(80,190){$7$}
\put(100,190){$2$}
\put(120,190){$5$}
\thicklines
\put(2,180){\line(1,0){80}}
\put(22,160){\line(1,0){40}}
\put(42,140){\line(1,0){20}}
\put(42,120){\line(1,0){60}}
\put(62,100){\line(1,0){20}}
\put(62,80){\line(1,0){60}}
\put(102,60){\line(1,0){20}}
\thinlines
\put(12,200){\line(0,-1){180}}
\put(32,200){\line(0,-1){180}}
\put(52,200){\line(0,-1){180}}
\put(72,200){\line(0,-1){180}}
\put(92,200){\line(0,-1){180}}
\put(112,200){\line(0,-1){180}}
\put(6,0){$\alpha_1$}
\put(26,0){$\alpha_2$}
\put(46,0){$\alpha_3$}
\put(66,0){$\alpha_4$}
\put(86,0){$\alpha_5$}
\put(106,0){$\alpha_6$}

\put(200,190){$5$}
\put(220,190){$3$}
\put(240,190){$1$}
\put(260,190){$7$}
\put(280,190){$4$}
\put(300,190){$2$}
\put(320,190){$6$}
\thicklines
\put(202,180){\line(1,0){60}}
\put(202,160){\line(1,0){120}}
\put(222,140){\line(1,0){40}}
\put(222,120){\line(1,0){60}}
\put(242,100){\line(1,0){20}}
\put(242,80){\line(1,0){40}}
\put(242,60){\line(1,0){60}}
\put(282,40){\line(1,0){40}}
\put(302,20){\line(1,0){20}}
\thinlines
\put(212,200){\line(0,-1){180}}
\put(232,200){\line(0,-1){180}}
\put(252,200){\line(0,-1){180}}
\put(272,200){\line(0,-1){180}}
\put(292,200){\line(0,-1){180}}
\put(312,200){\line(0,-1){180}}
\put(206,0){$\alpha_1$}
\put(226,0){$\alpha_2$}
\put(246,0){$\alpha_3$}
\put(266,0){$\alpha_4$}
\put(286,0){$\alpha_5$}
\put(306,0){$\alpha_6$}
\end{picture}
\]

\subsection{Necessity of the combinatorial conditions in Theorem~\ref{thm:main}}

It is possible to prove necessity by appealing to the geometric
description of the singularities along the maximal singular locus
found in \cite{Cortez, manivel2}; however we will give a simple, purely
combinatorial proof.

We will need the following two lemmas, the first of which is immediate:
\begin{Lemma}
\label{lemma:zero}
The vector $(\alpha_1,\ldots, \alpha_{n-1})\in {\mathbb Z}^{n-1}$ is a
solution to (\ref{eqn:final_system}) if and only if
$\sum_{r=i}^{j-1}\alpha_r =1$ for all $(i\leftrightarrow j)$ such that
$w(i\leftrightarrow j)$ covers $w$ in the Bruhat order.
\end{Lemma}

\begin{Lemma}
\label{lemma:seq_bars}
If there exists a solution $(\alpha_1,\ldots,\alpha_{n-1})\in {\mathbb
Z}^{n-1}$ to (\ref{eqn:final_system}), and $i<j$ with $w(i)<w(j)$, then
$\sum_{r=i}^{j-1}\alpha_r \geq 1$.  Equality holds if and only if
$w(i\leftrightarrow j)$ covers $w$.
\end{Lemma}
\begin{proof}
If $w(i\leftrightarrow j)$ covers $w$ then the claim holds by 
Lemma~\ref{lemma:zero}. Otherwise, it follows from the observation that
there are indices
\[i_0 = i <i_1<i_2 <\ldots < i_{t-1}<j=i_t\]
such that $w(i_s \leftrightarrow i_{s+1})$ covers $w$ for $0\leq s\leq t-1$.
\end{proof}

Now suppose that there is an embedding $i_1<i_2<i_3<i_4<i_5$ of a
$31524$ pattern with Bruhat restrictions 
$\{(1\leftrightarrow 5),(2\leftrightarrow 3)\}$.
Then by Lemma~\ref{lemma:seq_bars}, any solution would satisfy
\[\sum_{r=i_1}^{i_{3}-1}\alpha_r \geq 1, \sum_{r=i_2}^{i_4 -1}\alpha_r\geq 1,
\sum_{r=i_4}^{i_5 -1}\alpha_{r}\geq 1, \mbox{ and } \sum_{r=i_2}^{i_3 -1} \alpha_r =1.\]
Therefore, 
\begin{equation}
\label{eqn:askedforit}
\sum_{r=i_1}^{i_5 -1}\alpha_r =\sum_{r=i_1}^{i_3 -1}\alpha_r  +
\sum_{r=i_2}^{i_4 -1}\alpha_r + \sum_{r=i_4}^{i_5 -1}\alpha_r
-\sum_{r=i_2}^{i_3 -1}\alpha_r\geq 2.
\end{equation}
Thus (\ref{eqn:askedforit}) is a contradiction of
Lemma~\ref{lemma:zero} (or Lemma~\ref{lemma:seq_bars}) since $w(i_1
\leftrightarrow i_5)$ covers $w$. Therefore such an embedding cannot
exist.  A similar argument shows that there cannot exist an embedding
into $w$ of a $24153$ pattern with Bruhat restrictions
$\{(1\leftrightarrow 5),(3\leftrightarrow 4)\}$.

It remains to show that for each descent $d$ of $w$, 
$\lambda({\widetilde v}_{d}(w))$ has all of its inner corners on the same antidiagonal. For this purpose, we need:
\begin{Lemma}
\label{lemma:bij_inner}
Let $v$ be a Grassmannian permutation with descent at position
$d$. Then the transpositions $(i\leftrightarrow j)$ with $i\leq d<j$
such that $v(i\leftrightarrow j)$ covers $v$ are in bijection with the
inner corners of $\lambda(v)$.
Moreover, if $(i\leftrightarrow j)$ corresponds to an inner corner
of $\lambda(v)$ under this bijection, then the corresponding inner corner
distance equals $j-i-1$.
\end{Lemma}
\begin{proof}
In terms of the lattice path description of $\lambda(v)$ given on
page~2, an inner corner of $\lambda(v)$ occurs exactly when there is
an ``up step'' at time $a$, followed by a ``right step'' at time
$a+1$. In terms of $v$, this means $a$ and $a+1$ appear in positions
$i$ and $j$ satisfying the hypotheses.  Conversely, if $i\leq d < j$
and $v(i \leftrightarrow j)$ covers $v$, then $v(j)=v(i)+1$. The
claims then follow.
\end{proof}

The next lemma is clear from the definition of $v_{d}(w)$:

\begin{Lemma}
\label{lemma:same_both}
Let $d$ be a descent of $w$ and suppose $(i, j)$ is a pair $1\leq
i<j\leq n$ that indexes two entries of $w$ included in the subword
$v_d (w)$ of $w$.  Let $(i',j')$ be the corresponding indices in
${\widetilde v}_{d}(w)$. Then $w(i\leftrightarrow j)$ covers
$w$ if and only if ${\widetilde v}_{d}(w)(i'\leftrightarrow
j')$ covers ${\widetilde v}_{d}(w)$.
\end{Lemma}

Let $d$ be a descent of $w$ and suppose that 
\[i_1<i_2<\ldots <i_f =a <\ldots <i_s =d <i_{s+1}=d+1 <i_{s+2}<\ldots <i_{g}=b <\ldots < i_t\]
are the indices of the subword $v_{d}(w)$ of $w$, where
$w(a\leftrightarrow b)$ covers $w$. By Lemmas~\ref{lemma:seq_bars}
and~\ref{lemma:same_both} combined, any solution satisfies
\[1=\sum_{r=a}^{b-1}\alpha_r = (s-f) + (g-s-1) + \alpha_{d}=g-f-1+\alpha_{d}\]
Now, $g-f-1$ is the inner corner distance of the corresponding inner
corner of $\lambda({\widetilde v}_{d}(w))$ under the bijection of
Lemma~\ref{lemma:bij_inner}. Since $\alpha_{d}$ is fixed, $g-f-1$ is
independent of our choice of $a$ and $b$. Hence, all of the inner
corners of $\lambda({\widetilde v}_{d}(w))$ have the same inner corner
distance, and therefore they must all lie on the same antidiagonal.

\subsection{Sufficiency of the combinatorial conditions of Theorem~\ref{thm:main}}
Assume that the combinatorial conditions of Theorem~1 hold. 
We will show that in fact
\begin{equation}
\label{eqn:the_soln}
\alpha_r = \left\{
\begin{array}{cc}
1-{\mathfrak I}(\widetilde{v}_{r}(w)) & \mbox{ if $r$ is a descent} \\  
1 & \mbox{ otherwise.}
\end{array}
\right.
\end{equation}
for $1\leq r\leq n-1$ solves (\ref{eqn:final_system}).

It suffices to show that $\sum_{r=i}^{j-1} \alpha_r =1$ whenever
$w(i\leftrightarrow j)$ covers $w$. 
We prove this by induction on $j-i\geq 1$.

The base case $j-i=1$ of the induction holds by our definition of
$\alpha_r$, since in this case, $w$ does not have a descent at
position~$i$.

Now suppose that $j-i>1$. Let $k$ be chosen (if possible) so that
$i<k<j$ and $w(k)$ is minimal such that $w(k)>w(j)$.  Similarly, let
$\ell$ be chosen (if possible) so that $i<\ell<j$ and $w(\ell)$ is
maximal such that $w(\ell)<w(i)$. Notice that since
$w(i\leftrightarrow j)$ covers $w$, at least one of $k$ or $\ell$ must
exist. We now separately examine the possible cases:

First suppose $k$ exists but not $\ell$. Observe that $w$ has a
descent at position $j-1$, since, in fact $w(j)<w(m)$ for all $i\leq
m\leq j-1$. So we may consider the subword $v_{j-1}(w)$ of $w$.
Notice that this necessarily includes $w(i)$, $w(j-1)$, and $w(j)$.
By Lemma~\ref{lemma:same_both}, if $f$ and $g$ are indices between $i$
and $j-1$ in $w$ which correspond to successive entries of
$v_{j-1}(w)$, then $w(f\leftrightarrow g)$ covers $w$. So by induction,
\begin{equation}
\label{eqn:used_again}
\sum_{r=f}^{g-1}\alpha_{r} =1.
\end{equation}
Since by assumption, the inner corner distances of ${\widetilde v}_{j-1}(w)$ are all the same, by (\ref{eqn:the_soln}):
\[\sum_{r=i}^{j-1}\alpha_{r} = \sum_{r=i}^{j-2}\alpha_{r} + \alpha_{j-1} = {\mathfrak I}({\widetilde v}_{j-1}(w))+\alpha_{j-1} =1\]
as desired. A similar argument works in the case that $\ell$ exists but not
$k$, except that $v_{i}(w)$ must be used instead.

Next suppose that both $k$ and $\ell$ exist. First consider the 
situation where $k>\ell$. Then by construction, 
\[w(i\leftrightarrow k), \ w(\ell\leftrightarrow j),  
\mbox{ and $w(\ell\leftrightarrow k)$ each cover $w$}.\]
Therefore, by the induction hypothesis, we have
\[\sum_{r=i}^{k-1}\alpha_r = 1, \ \sum_{r=\ell}^{j-1} \alpha_r = 1, \mbox{ and }\sum_{r=\ell}^{k-1}\alpha_r = 1.\]
Hence,
\[\sum_{r=i}^{j-1}\alpha_r = \sum_{r=i}^{k-1}\alpha_r + \sum_{r=\ell-1}^{j-1}\alpha_r - \sum_{r=\ell}^{k-1}\alpha_r=1\]
as desired.

Finally, we have the case where $k<\ell$. Observe that the values of
$w$ between $k$ and $\ell$ must consist of numbers larger than $w(k)$
followed by numbers smaller than $w(\ell)$, since otherwise it is easy
to see that there must exist a $\{(1\leftrightarrow 5),
(2\leftrightarrow 3)\}$-restricted embedding of $31524$ or a
$\{(1\leftrightarrow 5), (3\leftrightarrow 4)\}$-restricted embedding
of $24153$, contradicting the assumptions.  Similarly, the values of
$w$ between $i$ and $k$ are necessarily smaller than $w(i)$.

Let $q$ be the last index $k\leq q<\ell$ such that $w(q)\geq w(k)$;
hence $w$ has a descent at~$q$.  Consider the subword $v_{q}(w)$ of
$w$ and observe that $w(i)$ and $w(j)$ are in $v_{q}(w)$, as,
otherwise, we would find a bad $31524$ or $24153$ pattern.  We are now
ready to employ a similar argument as above.  By
Lemma~\ref{lemma:same_both}, if $f$ and $g$ are indices of $w$, with
either both $f$ and $g$ in the interval $[i,q]$ or both in the
interval $[q+1,j]$, and $f$ and $g$ correspond to consecutive entries
of $v_{q}(w)$, then $w(f\leftrightarrow g)$ covers $w$; now the
induction hypothesis implies (\ref{eqn:used_again}) as before.
Therefore, by (\ref{eqn:the_soln}) and our assumptions about
$\lambda({\widetilde v}_{q}(w))$, we have
\[\sum_{r=i}^{j-1}\alpha_r = {\mathfrak I}({\widetilde v}_{q}(w))+\alpha_q =1\]
as required.

Theorem~\ref{thm:main} follows immediately
from the discussion above.

\subsection{Conclusion of the proof of Theorem~\ref{thm:second_main}}   
In order to complete the above arguments to prove
Theorem~\ref{thm:second_main}, we need two facts about the Borel-Weil
construction; see, for example,~\cite[Section~1.4]{brion:book} and the
references therein.  First, we note that, if
$\mathcal{L}_{\Lambda_{n-r}}$
denotes the line bundle associated to the fundamental
weight $\Lambda_{n-r}$ by the Borel-Weil construction, then
$\mathrm{div}(\mathcal{L}_{\Lambda_{n-r}})=\left[X_{(r\leftrightarrow
r+1)}\right]\in {\rm CaCl}({\rm Flags}(C^n))$; therefore,
\[\mathrm{div}(\mathcal{L}_{\Lambda_{n-r}}\big|_{X_w})=\sum_{\plainstack{a\leq r<b}{\ell(w(a\leftrightarrow b))=\ell(w)+1}} \left[ X_{w(a\leftrightarrow
b)}\right]\in {\rm CaCl}(X_w).\] (The line bundle
$\mathcal{L}_{\lambda_{n-r}}$ can be concretely constructed using the
isomorphism $\mathcal{L}_{\Lambda_{n-r}}\cong \bigwedge^{n-r}
\mathcal{Q}_r$, , where $\mathcal{Q}_r$ is the tautological quotient
bundle whose fiber at a flag $F_{\bullet}=(\langle 0\rangle \subseteq
F_1\subseteq\ldots\subseteq F_n={\mathbb C}^n)$ is $\mathbb{C}^n/F_r$.)
Secondly, addition of weights corresponds to tensor product of line
bundles, so that, for any weights $\alpha$ and $\beta$, the line
bundle
$\mathcal{L}_{\alpha+\beta}=\mathcal{L}_\alpha\otimes\mathcal{L}_\beta$.

We have shown that, when $X_w$ is Gorenstein,
$$\sum_{r=1}^{n-1}
\alpha_r \ \mathrm{div}(\mathcal{L}_{\Lambda_{n-r}}\big|_{X_w}) =
\sum_{\plainstack{v\succ w}{\ell(v)=\ell(w)+1}}
\left[X_v\right]=\mathrm{div}(\mathcal{I}(\partial X_w)).$$ Therefore, we have
that $\mathcal{I}(\partial X_w)\cong \mathcal{L}_\alpha\big|_{X_w}$, where
$\alpha=\sum_{r=1}^{n-1}-\alpha_r\ \Lambda_{n-r}$.  Since
$\rho=\sum_{r=1}^{n-1}\Lambda_r$, and we have set
$\widetilde{\alpha}_r=-1-\alpha_r$ in (\ref{eqn:the_soln'}),
this proves Theorem~\ref{thm:second_main}.  \qed

\subsection{Proof of Corollary~\ref{cor:maxsingloc}}

We prove Corollary~\ref{cor:maxsingloc} by comparing
Theorem~\ref{thm:main} with a description of the generic singularities
of a Schubert variety given in~\cite{Cortez, manivel2}.~\footnote{Note that
our notation differs from the notation in these papers by right multiplication of a permutation
$w$ by $w_0$.}

Suppose a Schubert variety $X_w$ is not Gorenstein along its maximal
singular locus.  By the local definition of Gorensteinness given in
Section~2.1, it is not Gorenstein.  To prove the other direction,
suppose $X_w$ is not Gorenstein.  Then $w$ contains one of the two
forbidden patterns, or violates the inner corner condition.  If $w$
contains a forbiddern pattern, then, in the language of
Cortez~\cite{Cortez}, $w$ has a configuration II with $r=0$ and
$s+t\geq 1$, and therefore has a generic singularity whose
neighborhood is isomorphic to the product of $\mathbb{C}^k$ for some
$k$ and the variety of $(s+t+2)\times 2$ matrices of rank at most 1.
It is well known that the variety of $p\times q$ matrices of rank at
most 1 is Gorenstein if and only if $p=q$; see for
example~\cite[Thm. 7.3.6]{Bruns-Herzog}; this shows that $X_w$ is not
Gorenstein at a generic singularity.  If $w$ violates the inner corner
condition, then $w$ has a configuration I with $s\neq t$, yielding a
corresponding generic singularity, which, as it a neighborhood
isomorphic to the product of $\mathbb{C}^k$ for some $k$ and the
variety of $s\times t$ matrices of rank at most 1, is not Gorenstein.

\section{Remarks and Applications}

\subsection{Extension to partial flag varieties}
\label{subsection:partial}
More generally, let ${\rm Flags}(i_1<i_2<\ldots<i_k,{\mathbb C}^n)$
denote the variety of partial flags $F_{\bullet}:\langle 0\rangle
\subseteq F_{i_1} \subseteq F_{i_2}\subseteq \ldots\subseteq
F_{i_{k}}\subseteq {\mathbb C}^{n}$ in ${\mathbb C}^{n}$ where here
${\rm dim}(F_{i_{k}})=i_{k}$.  By convention let $i_0=0$ and
$i_{k+1}=n$.  Now let $S=S_{i_1-i_0}\times
S_{i_2-i_1}\times\cdots\times S_{i_{k+1}-i_k}\subseteq S_n$ denote the
Young subgroup where the $S_{i_j-i_{j-1}}$ factor is generated by the
simple reflections $s_{i_{j-1}+1},\ldots,s_{i_j -1}$ for
all $j$ such that $1\leq j\leq k$.  The Schubert varieties of ${\rm
Flags}(i_1<i_2<\ldots<i_k,{\mathbb C}^n)$ are indexed by cosets
of~$S$. The natural ``forgetting subspaces'' projection $\pi:{\rm Flags}({\mathbb C}^n)\twoheadrightarrow {\rm Flags}(i_1<i_2<\ldots<i_k,{\mathbb C}^n)$ 
is a smooth fiber bundle.  It follows that a
Schubert variety $X_{wS}$ in ${\rm Flags}(i_1<i_2<\ldots<i_k,{\mathbb
C}^n)$ indexed by a coset $wS$ is Gorenstein if and only if the
Schubert variety $X_{\widetilde w}=\pi^{-1}(X_{wS})$ in ${\rm
Flags}({\mathbb C}^n)$ is Gorenstein, where ${\widetilde w}$ is the
minimal length element of $wS$. In particular, our main result implies
the Grassmannian case as presented in the introduction.

\subsection{Uniqueness of (\ref{eqn:the_soln})} 
It is worthwhile to note that the induction in Section~2.4 implies that
(\ref{eqn:the_soln}) is a solution to (\ref{eqn:final_system}) if and only if
$X_w$ is Gorenstein. Moreover, this solution is essentially
unique. The only exception to uniqueness arises for those $r$ where 
\[\left[X_{(r\leftrightarrow r+1)}\right]\left[X_w\right] 
=\sum_{\plainstack{a\leq r<b}{\ell(w(a\leftrightarrow b))=\ell(w)+1}} \left[
X_{w(a\leftrightarrow b)}\right]=0\] because the sum on the right hand
side is vacuous.  In these cases, we can arbitrarily assign a value to
$\alpha_r$ in order to arrive at a solution.  (This is also apparent
from the bar diagrams of Section~2.2, as in these cases no horizontal
bars cross the $r^{th}$ vertical bar.)  Consequently, the expression
for $\omega_{X_w}$ given in Theorem~\ref{thm:second_main} is unique,
up to tensoring by bundles which are trivial when restricted to
$X_{w}$. Furthermore:

\noindent
\emph{${\mathbb Q}$-Gorensteinness:} A variety is said to be {\bf
${\mathbb Q}$-Gorenstein} if it is Cohen-Macaulay and some multiple of
the canonical divisor is Cartier.  Consequently, a Schubert variety
$X_w$ is ${\mathbb Q}$-Gorenstein if (\ref{eqn:final_system}) has a
rational solution.  However, since if any solution exists, an integral
solution exists, Gorensteinness and $\mathbb{Q}$-Gorensteinness are
equivalent.  This will not hold in general for flag varieties of other
Lie types.

\noindent
\emph{Computational efficiency:} 
In order to check if a permutation $w$ corresponds to a Gorenstein
Schubert variety, it is typically more computationally efficient solve for
(\ref{eqn:final_system}) than to use Theorem~\ref{thm:main}. In particular, it is enough to check if (\ref{eqn:the_soln}) works.

\subsection{Is it pattern avoidance?}
	In view of~\cite{LS90}, 
it is natural to wonder if it is possible to reformulate 
Theorem~1 in terms of ``classical pattern avoidance'', that is, if 
there is a finite list of permutations 
$w_1, w_2, \ldots, w_n$ such that $X_{w}$ is Gorenstein if and only if 
$w$ pattern avoids these permutations.

	In fact, this is already impossible for Grassmannian
permutations. For example, we know $1346\ \mid \ 25\in S_6$ does not
correspond to a Gorenstein Schubert variety. But $w' =
\underline{1}2\underline{5}\underline{6}\underline{9}\ \mid \
\underline{3}4\underline{7}8\in S_9$ does. Note that $w'$ contains $w$
as a subpattern, so if a classical pattern avoidance permutation
reformulation of Theorem~\ref{thm:main} existed, it would imply that
$X_{w'}$ is not Gorenstein, which is not true.

\subsection{A characterization of Fano Schubert varieties}

        A Gorenstein algebraic variety is {\bf Fano} if its
anticanonical divisor is ample.  It follows from
Theorem~\ref{thm:second_main} that a Gorenstein Schubert variety
$X_{w}$ in ${\rm Flags}({\mathbb C}^n)$ is Fano if and only if all of
the inner corner distances of $w$ are at most~1.  This appears to give
new examples of Fano varieties.  It seems to have been previously
unknown whether or not all smooth Schubert varieties of the flag
variety are Fano.  By the above remark, it is easy to find examples of
Schubert varieties that are smooth but not Fano, in contrast to the
case for Grassmannians, for which all smooth Schubert varieties are
Fano.

\subsection{Matrix Schubert varieties and ladder determinantal varieties}

Let $v\in S_n$ be a permutation, and $Y_v$ the associated matrix
Schubert variety; this was defined in \cite{Fulton:Duke92} as the
closure in $\mathbb{C}^{n^2}$, considered as the space of $n \times n$
matrices, of $p^{-1}(X_v)$, where $p: \mathrm{GL}_{n}({\mathbb C})
\rightarrow
\mathrm{GL}_n(\mathbb{C})/B=\mathrm{Flags}(\mathbb{C}^{n})$ is the
quotient map.  Now let $w=v \times {\rm id}\in S_n \times S_n \subseteq
S_{2n}$ be the permutation agreeing with $v$ on $1,\ldots, n$ and
fixing $n+1,\ldots, 2n$.  The intersection of $X_w$ with the opposite
big cell of flags intersecting the canonical reference flag (whose
$i$-th vector space is $\langle e_1,\ldots,e_i\rangle$) generically is
then isomorphic to $Y_v \times \mathbb{C}^{n^2-n}$.  Every singularity
of $X_w$ is represented in this opposite big cell, so $Y_v$ is
Gorenstein if and only if $X_w$ is.  Identifying ladder determinantal
varieties with the appropriate matrix Schubert varieties allows us to
recover the characterizations of Gorenstein ladder determinantal
varieties found in \cite{Conca} and \cite{GonMiller}.

\subsection{Theorem~\ref{thm:second_main} and cohomology of line bundles 
on Gorenstein Schubert varieties}

Theorem~\ref{thm:second_main} can be applied to obtain information
about the sheaf cohomology groups ${\rm H}^{i}(X_w, {\mathcal
L}_{\alpha}\big|_{X_w}\big)$ of the line bundle ${\mathcal
L}_{\alpha}\big|_{X_w}$ on a Gorenstein Schubert variety $X_w$. The
groups are classically known in the case $X_{\rm id}\cong {\rm
Flags}({\mathbb C}^n)$ and $\alpha\in {\rm Hom}(T,{\mathbb C}^n)$ is
arbitrary (the classical Borel-Weil-Bott theorem~\cite{Bott}), and for
arbitrary $w\in S_n$ when $\alpha$ is dominant~\cite{Demazure};
see, for example,~\cite{Jantzen}.
It is an open problem to compute these groups
in most of the remaining cases;
see~\cite{transform} for some recent progress on this problem.

Serre duality (see, for example,~\cite[III.7]{Hartshorne}) 
states that, for any projective, equidimensional,
$d$-dimensional, Cohen-Macaulay scheme $X$, and any coherent sheaf
$\mathcal{F}$ on $X$, we have
\[H^i(X,\mathcal{F}) \cong\mathrm{Ext}^{d-i}(\mathcal{F},\omega_X)^* .\]
Let $\alpha$ be the (non-dominant) weight
defined in Theorem~\ref{thm:second_main}, and $\beta$ 
any weight.  Then: 
\begin{eqnarray*}
H^i\big(X_w,\mathcal{L}_{\alpha-\beta}\big|_{X_w}\big) & \cong & 
\mathrm{Ext}^{n-\ell(w)-i}(\mathcal{L}_{\alpha-\beta}\big|_{X_w},\omega_X)^* \\
& = & \mathrm{Ext}^{n-\ell(w)-i}(\mathcal{L}_{\alpha-\beta}\big|_{X_w},
\mathcal{L}_{\alpha}\big|_{X_w})^* \\
& \cong& \mathrm{Ext}^{n-\ell(w)-i}(\mathcal{O}_{X_w},\mathcal{L}_{\beta}\big|_{X_w})^* \\
& \cong & H^{n-\ell(w)-i}\big(X_w,\mathcal{L}_{\beta}\big|_{X_w}\big)^*.
\end{eqnarray*}
When $\beta$ is dominant, this relates the cohomology groups $H^i\big(X_w,\mathcal{L}_{\alpha-\beta}\big|_{X_w}\big)$
to the cohomology groups known by Demazure's theorem. For example, 
it follows that, when $\beta$ is dominant, $H^i\big(X_w,\mathcal{L}_{\alpha-\beta}\big|_{X_w}\big)\cong 0$ for 
$i\neq n-\ell(w)$.


\end{document}